\documentclass[11pt]{sat}

\newif\ifsattoc\sattoctrue
\newread\testfl\immediate\openin\testfl=\jobname.tic
\ifeof\testfl\sattocfalse\fi\closein\testfl

\title{On the Power of Function Values for the Approximation\vskip0pt
Problem in Various Settings}

\def\shorttitle{The Power of Function Values}

\author{Erich Novak and Henryk Wo\'zniakowski}
\def\shortauthor{E. Novak, H. Wo\'zniakowski}

\def\versiondate{June 2, 2011}

\def\abstracttext{This
is an expository paper on
approximating functions from general
Hilbert or Banach spaces
in the worst case, average case and randomized settings with error
measured in the~$L_p$ sense.
We define the~\emph{power function} as
the ratio between the best rate of convergence of algorithms
that use function values over the best rate of convergence
of algorithms that use arbitrary linear functionals
for a worst possible Hilbert or Banach space for which the
problem of approximating functions is well defined.
Obviously, the power function takes values at most one.
If these values are one or close to one
than the power of function values is the same or almost the same
as the power of arbitrary linear functionals.
We summarize and supply a few new estimates on the power function.
We also indicate eight open problems related to the power function
since this function has not yet been studied in many cases.
We believe that the open problems will be of interest
to a general audience of mathematicians.}

\def\MSCnumbers{41A25, 41A46, 65Y20}

\def\keywords{degree of approximation, widths and sampling numbers,
complexity of numerical algorithms}

\usepackage{latexsym,amsfonts,amsmath,epsfig,tabularx,amsthm}

\theoremstyle{plain}

\newtheorem{thm}{Theorem}

\theoremstyle{definition}

\newtheorem{rem}{Remark}

\newtheorem{example}{Example}
\newtheorem{OP}{Open Problem}
\newtheorem{defi}{Definition}

\newcommand\expect{\mathbb{E}}
\newcommand\reals{\mathbb{R}}
\newcommand\R{{\mathbb{R}}}

\newcommand\nat{\mathbb{N}}
\newcommand\il{\left<}
\newcommand\ir{\right>}
\newcommand\eps{\varepsilon}

\newcommand\lall{\Lambda^{{\rm all}}}
\newcommand\lstd{\Lambda^{{\rm std}}}

\newcommand{\e}{\varepsilon}
\newcommand{\enwall}{e^{\rm all-wor}_n}
\newcommand{\enwstd}{e^{\rm std-wor}_n}
\newcommand{\enaall}{e^{\rm all-avg}_n}
\newcommand{\enastd}{e^{\rm std-avg}_n}
\newcommand{\enrall}{e^{\rm all-ran}_n}
\newcommand{\enrstd}{e^{\rm std-ran}_n}

\newcommand{\infp}{\operatornamewithlimits{inf\phantom{p}}}

\newcommand{\wt}{\widetilde}

\def\startpagenumber{1}
\def\volumenumber{6} 
\def\year{2011}

\newcommand{\beginddoc}{
\begin{document}
\maketitle
\begin{abstract}
\abstracttext
\vskip1pt MSC: \MSCnumbers
\ifx\keywords\empty\else\vskip1pt Keywords: \keywords\fi
\end{abstract}
\insert\footins{\scriptsize
\medskip
\baselineskip 8pt
\leftline{Surveys in Approximation Theory}
\leftline{Volume \volumenumber, \year.
pp.~\thepage--\pageref{endpage}.}
\leftline{\copyright\ \year\ Surveys in Approximation Theory.}
\leftline{ISSN 1555-578X}
\leftline{All rights of reproduction in any form reserved.}
\smallskip
\par\allowbreak}
\ifsattoc\else\tableofcontents\fi}
\renewcommand\rightmark{\ifodd\thepage{\it \hfill\shorttitle\hfill}\else
{\it \hfill\shortauthor\hfill}\fi}
\markboth{{\it \shortauthor}}{{\it \shorttitle}}
\markright{{\it \shorttitle}}
\def\endddoc{\label{endpage}\end{document}}
\date{{\small \versiondate}}
\setlength\oddsidemargin{0pc}
\setlength\evensidemargin{0pc}
\setlength\topmargin{0in}
\setlength\textwidth{6.5in}
\setlength\textheight{8.6in}
\beginddoc
\ifsattoc
\bigskip
\def\toczer{0}\def\tochalf{.5}\def\tocone{1}
\def\tocindent{0}
\def\ection{section}\def\ubsection{subsection}
\def\numberline#1{\hskip\tocindent truecm{} #1\hskip1em}
\newread\testfl
\def\inputifthere#1{\immediate\openin\testfl=#1
    \ifeof\testfl\message{(#1 does not yet exist)}
    \else\input#1\fi\closein\testfl}
\countdef\counter=255
\def\diamondleaders{\global\advance\counter by 1
  \ifodd\counter \kern-10pt \fi
  \leaders\hbox to 15pt{\ifodd\counter \kern13pt \else\kern3pt \fi
  \hss.\hss}\hfill}
\newdimen\lextent
\newtoks\writestuff
\medskip
\begingroup
\small
\def\contentsline#1#2#3{
\def\argu{#1}
\ifx\argu\ection\let\tocindent\toczer\else
\ifx\argu\ubsection\let\tocindent\tochalf\else\let\tocindent\tocone\fi\fi
\setbox1=\hbox{#2}\ifnum\wd1>\lextent\lextent\wd1\fi}
\lextent0pt\inputifthere{\jobname.tic}\advance\lextent by 2em\relax
\def\contentsline#1#2#3{
\def\argu{#1}
\ifx\argu\ection\let\tocindent\toczer\else
\ifx\argu\ubsection\let\tocindent\tochalf\else\let\tocindent\tocone\fi\fi
\writestuff={#2}
\centerline{\hbox to \lextent{\rm\the\writestuff%
\ifx\empty#3\else\diamondleaders{}
\hfil\hbox to 2 em\fi{\hss#3}}}}
\input \jobname.tic\endgroup \newpage
\fi



\section{Introduction}  \label{s1}

This is an expository paper on the
problem of approximating functions from general
Hilbert or Banach spaces, which
has been thoroughly studied in many books and papers.
This problem has many variants depending on how we measure the error
of such approximations (algorithms).
A popular choice is to take the norm of an
$L_p$ space and all values of $p\in[1,\infty]$ have been considered.
Furthermore, the error of algorithms
can be defined in the worst case, average case
or randomized setting.
For the worst and average case settings, we consider
deterministic algorithms. The worst case error is defined
as the maximal error over the unit ball of a given space whereas the
average case error is defined as the average error over the whole
space with respect to a given measure. The usual choice
is a zero mean Gaussian measure. For the randomized setting
we consider randomized algorithms and the error is defined as the
maximal expected error over the unit ball of a given space.
Here, the expected error is given with respect
to a probability distribution of randomized elements.

We approximate functions $f$ by algorithms that use information about $f$
given by  finitely many
functionals of $f$. Information is called~\emph{linear}
if we can choose arbitrary linear functionals, and it is
called~\emph{standard} if only function values
may be used. Clearly, linear information is at least as powerful as
standard information. For many applications, only standard information
is available. But even in this case, it is a good idea to study linear
information and learn how difficult is the function approximation
problem. For example, if we can prove that even for linear information
the problem is too difficult then, obviously, the same also holds for
standard information. On the other hand, all positive results for
linear information do not have to hold for standard information.

The main question addressed in this expository paper
is the study of the power of standard information or equivalently
the power of function values.
We want to know how much we lose if function
values are used instead of linear information. Or more optimistically,
we ask when the power of standard information is the same
or nearly the same as the power of linear information.
Such questions have been addressed in a number of papers and we will
refer to them in the course of this paper. It has been usually done
for specific spaces and only a few papers addressed these questions
for some classes of spaces.

Our approach is a little more general and we want to verify the power
of function values/standard information for all Hilbert or Banach
spaces for which the problem of function approximation is well
defined. More precisely, we define
the~\emph{power} function\footnote{We needed to find a good one-letter
name for the power function. Since in English and in Polish this
would indicate the letter ``$p$'' which is already used as the parameter
of the $L_p$ space, we turn to German and use the word ``Leistung''.
That is why the letter $\ell$ denotes the power function.}
$$
\ell^{\,\rm sett-x}\,:\,(0,\infty)\times [1,\infty]\to [0,1].
$$
Here ${\rm sett}\in\{{\rm wor},{\rm ran},{\rm avg}\}$ denotes the
setting we use for the error definition. Hence, ${\rm wor}$ stands for
the worst case setting, ${\rm ran}$ for the randomized setting, and
${\rm avg}$ for the average case setting.
The second superscript ${\rm x}\in\{H,B\}$ tells us if we consider only
Hilbert spaces (${\rm x}=H$) or if we allow all Banach spaces
(${\rm x}=B$).

We now explain the meaning of the value
$$
\ell^{\,\rm sett-x}(r,p).
$$
The first argument $r$ means that
the $n$th minimal error (formally defined in Definition 1)
behaves like $n^{-r}$
if we use linear information.
Since $r>0$, we consider Hilbert or Banach spaces which admit
convergence, and furthermore they admit a
polynomial rate of convergence of the minimal errors.
The second argument $p$ denotes the use of the norm of $L_p$.
The value
$\ell^{\,\rm sett-x}(r,p)$ is defined as $r^{-1}$ times the best rate of
convergence we obtain using only function values
for a worst possible choice of a Hilbert or Banach space.
That is why $\ell^{\,\rm sett-x}(r,p)\le1$, and the larger
$\ell^{\,\rm sett-x}(r,p)$ the better.
Hence, if we have
$$
\ell^{\,\rm sett-x}(r,p)=1
$$
then the power of standard information
is the same as the power of linear information.
We will see later that this
does happen in some cases.
Then
standard information yields the same rate
of convergence as linear information for
the embeddings $I: F \to L_p$ for \emph{all} Hilbert (if $x=H$)
or Banach (if $x=B$) spaces
without the need of a case to case study for each $F$.
This holds in the randomized setting for Hilbert spaces
with $p=2$, see Theorem~\ref{thm6}, and in the average case setting
for Banach spaces equipped with zero mean Gaussian measures and
$p=2$, see Theorem~\ref{thm8}. It is open if
$\ell^{\, \rm sett-x}(r,p)=1$ may happen in the worst case setting,
see Open Problem 1.

On the other hand,
if we have
$$
\ell^{\,\rm sett-x}(r,p)=0
$$
then the power of standard information is zero as
compared to the power of linear information. Finally, if we have
$$
\ell^{\,\rm sett-x}(r,p)\in(0,1)
$$
then we know qualitatively how much we may lose by using function values.

The concept of the power function seems to be new. For many values of
$p$, especially when $p\not=2$, this function has not yet been
studied. This is especially the case for the randomized and average
case settings. That is why we indicate eight open problems related
to the power function with the hope that many mathematicians
will be interested in solving them and advancing our knowledge
about the power of function values.

In this paper, we tried to summarize and supply a few new
estimates on the power function. We now briefly indicate
a few results presented in the paper.

In the worst case setting for the Hilbert case and $p=2$,
we conclude from~\cite{HNV08,KWW09} that
\begin{eqnarray*}
\ell^{\,\rm wor-H}(r,2)&=&0\ \ \ \ \ \ \ \ \ \ \ \ \
\ \ \ \ \ \mbox{for all}\ \ \
r\in(0,\tfrac12],\\
\ell^{\,\rm wor-H}(r,2)&\in&\left[\frac{2r}{2r+1},1\right]
\ \ \ \ \, \mbox{for all}\ \ \
r\in(\tfrac12,\infty).
\end{eqnarray*}
Hence, the power of function values is zero for $r\le1/2$,
and almost the same as the power of linear information for large $r$.
One of the main open problem is to verify whether
$\ell^{\,\rm wor-H}(r,2)=1$ for all $r>1/2$.

Staying with the worst case and Hilbert spaces but with $p\not=2$,
we conclude from~\cite{Ta10} that
$$
\ell^{\,\rm wor-H}(r,p)=0\ \ \ \mbox{for all}\ \ \
r\in\left(0,\min(\tfrac1p,\tfrac12)\right].
$$
For $r>\min(1/p,1/2)$, we do not know anything about
the values of $\ell^{\,\rm
wor-H}(r,p)$ except the case $p=\infty$ for which we know
from~\cite{No88} that
$$
\ell^{\,\rm wor-H/B}(r,\infty)\ge1-\frac1r.
$$
By H/B we mean that we obtain this result for both Hilbert and
Banach spaces.
Again for large~$r$, the power of standard information is almost the
same as the power of linear information.

For the worst case and the Banach case, we have
\begin{eqnarray*}
\ell^{\,\rm wor-B}(r,p)&=&0\ \ \ \ \ \ \ \ \
\ \ \ \ \ \ \ \qquad \quad
\mbox{for all}\ \ \
r\in(0,1]\ \ \mbox{and}\ \ p\in [1, 2]    ,\\
\ell^{\,\rm wor-B}(r,p)&=&0\ \ \ \ \ \ \ \ \
\ \ \ \ \ \ \ \qquad \quad
\mbox{for all}\ \ \
r\in(0, \tfrac{1}{2} + \tfrac{1}{p}] \ \ \mbox{and}\ \ p\in(2, \infty )
     ,\\
\ell^{\,\rm wor-B}(r,p)&\le& 1-\frac1r\left(1-\frac1p\right)
\ \ \ \ \ \ \, \mbox{for all}\ \ \ r>1\ \ \mbox{and}\ \ p\in[1,2],\\
\ell^{\,\rm wor-B}(r,p)&\le& 1-\frac1{2r}\ \ \ \ \ \ \ \ \ \ \ \ \
\ \ \ \ \ \mbox{for all}\ \ \ r>1\ \ \mbox{and}\ \ p\in[2,\infty), \\
1- \frac{1}{r} \le \ell^{\,\rm wor-B}(r,\infty)&\le& 1-\frac1{2r}\ \ \
\ \
\ \ \ \ \ \ \  \ \  \ \ \  \ \mbox{for all}\ \ \ r>1 ,
\end{eqnarray*}
       see Theorem~\ref{thm5} in Section~\ref{s2.2}.

Even though we do not know much about the power function in this case,
we can conclude that the Hilbert and Banach cases are different since
$$
\ell^{\,\rm wor-B}(r,2)<\ell^{\,\rm wor-H}(r,2)
\ \ \ \mbox{for all}\ \ \ r\in(\tfrac12,\infty).
$$

Surprisingly enough, for the randomized setting with the Hilbert
case and for the average case setting with the Hilbert or Banach case
we have complete knowledge about the power function for
$p=2$ due to~\cite{WW07} and~\cite{HWW08}.
More precisely, we know that
$$
\ell^{\,\rm ran-H}(r,2)=\ell^{\,\rm avg-H/B}(r,2)=1
\ \ \ \mbox{for all}\ \ \ r>0.
$$
More estimates of the power function can be found in the
    subsequent sections.

\section{Worst case setting}  \label{s2}

Let $F$ be a Hilbert or Banach space
of functions, defined on a set $\Omega$,
such that the linear functionals
$f \mapsto f(x)$ are continuous for all $x\in
\Omega$.
We assume that
$F \subset L_p$ and that the embedding
$I: F \to L_p$ is continuous\footnote{We do not specify
$\Omega$ or the underlying measure of $L_p$ since they
can be arbitrary.}, where $I(f)=f$.
We write $H$ instead of $F$ if $F$ is a Hilbert space.

Let $(c_n)$ be a sequence of nonnegative numbers.
Assume first that $(c_n)$ converges to zero.
We define its (polynomial) rate of convergence $r(c_n)$ by
$$
r(c_n) = \sup \{\, \beta \ge 0 \mid\lim_{n \to \infty} c_n n^\beta = 0 \}.
$$
If $(c_n)$ is not convergent to zero, we set $r(c_n)=0$.
Then $r(c_n)$ is well defined for all nonnegative
sequences $(c_n)$. For example,
the rate of convergence of $n^{-\alpha}$ is $\max(0,\alpha)$.

We approximate functions from $F$ using finitely many arbitrary linear
functionals $L\in F^*$ or function values $f(x)$ for some $x\in
\Omega$.
We define the error
of such approximations by taking the worst case setting with respect to
the $L_p$ norm. The norm of $L_p$ is denoted by $\|\cdot\|_p$.

We define two classes $\lall$ and $\lstd$ of information evaluations.
We have $\lstd\subseteq\lall=F^*$ and $\lstd$ consists
of linear functionals
of the form $L_x(f)=f(x)$
for all $f\in F$, where $x\in\Omega$.  We approximate functions from
$F$ by algorithms $A_n:F\to L_p$ given by
$$
A_n(f)=\phi_n(L_1(f),L_2(f),\dots,L_n(f)),
$$
where $n$ is a nonnegative integer, $\phi_n:\reals^n\to L_p$ is an
arbitrary mapping, and $L_j\in \Lambda$, where
$\Lambda \in \{\lall, \lstd\}$.
The choice of $L_j$ can be adaptive, that is,
$L_j(\cdot)=L_j(\cdot;L_1(f),L_2(f),\dots,L_{j-1}(f))$
may depend on the
already computed values $L_1(f),L_2(f),\dots,L_{j-1}(f)$.
For $n=0$, $A_n(f)$ equals
some fixed element of the space $L_p$.
More details can be found in e.g.,~\cite{NW08,TWW88}.

Hence, we consider algorithms that use $n$ linear functionals either
from the class $\lstd$ or from the class $\lall$.
We define the minimal errors as follows.

\begin{defi}
For $n=0$ and $n\in \nat:=\{1,2,\dots \}$, let
$$
\enwall(F, L_p) := \infp_{A_n \ {\rm with}\ L_j\in \lall}
\  \sup_{\Vert f \Vert_F \le 1}
\big\Vert f - A_n(f)\big \Vert_p
$$
and
$$
\enwstd(F, L_p) := \infp_{A_n\ {\rm with}\ L_j\in \lstd}
\ \sup_{\Vert f \Vert_F \le 1}
\big\Vert f - A_n(f)\big\Vert_p  .
$$
\end{defi}

For $n=0$, it is easy to see that the best
algorithm is $A_0(f)=0$ and we obtain
$$
e^{\rm all-wor}_0(F,L_p)= e^{\rm std-wor}_0(F, L_p)=
\sup_{\Vert f \Vert_F \le 1} \Vert f \Vert_p =
\sup_{\Vert f \Vert_F \le 1} \Vert I(f) \Vert_p =
\Vert I \Vert.
$$
This is the initial error that can be achieved without
computing any linear functional on the functions $f$.
Clearly,
$$
\enwall(F, L_p ) \le \enwstd(F, L_p)\ \ \ \mbox{for all} \ \ \  n\in\nat.
$$
The sequences $\left(\enwall(F, L_p)\right)$
and $\left(\enwstd(F, L_p)\right)$
are both non-increasing but not necessarily convergent to zero.

We want to compare the rates of convergence
$$
r^{\rm all-wor}(F,L_p)  :=
r\left(\enwall(F,L_p)\right) \qquad
\hbox{and}  \qquad
r^{\rm std-wor}(F,L_p) := r\left(\enwstd(F,L_p)\right) .
$$
In particular, we would like to know if it is
possible that the sequence $\left(\enwall(F,L_p)\right)$
converges to zero much faster than the sequence
$\left(\enwstd(F,L_p)\right)$.
In many cases it is much easier to analyze
the sequence $(\enwall(F,L_p))_{n \in \nat}$.
It is then natural to ask what can be said about
the sequence $(\enwstd(F,L_p))_{n\in\nat}$.

The main question addressed in this paper is to find or estimate
the~\emph{power} function defined as
$\ell^{\,\rm wor-x}:(0,\infty)\times[1,\infty]\to[0,1]$ by
$$
\ell^{\,\rm wor-x}(r,p)
:=\inf_{F:\, r^{\rm all-wor}(F,L_p)=r
}\frac{r^{\rm std-wor}(F,L_p)}r,
$$
where ${\rm x}\in\{H,B\}$ and indicates that
the infimum is taken over all Hilbert spaces $({\rm x}=H)$ or over all
Banach spaces $({\rm x}=B)$ continuously embedded in $L_p$
for which function values are continuous linear functionals
and the rate of convergence is $r$ when we use
arbitrary linear functionals.

It is easy to show, and it will be shown later,
that the set of spaces $F$ for which $r^{\rm
  all-wor}(F,L_p)=r$ is not empty and therefore
$\ell^{\,\rm wor-x}$ is well defined.
Obviously, $\ell^{\,\rm wor-x}(r,p)\in[0,1]$, as already claimed.
The power function  $\ell^{\,\rm wor-x}$ measures the ratio
between the best rates of convergence of approximations based
on function values over those based on arbitrary linear
functionals for a worst
possible Hilbert or Banach space.

We briefly comment on why we take the infimum over $F$ in the definition
of the power function. For some specific spaces $F$,
standard information is as powerful as linear
information\footnote{This holds with $r=\infty$ for all finite
dimensional spaces $F$. This also holds for some infinite
dimensional Hilbert spaces $F$. For example, take $F$ as the space
of piecewise constant functions over, say, $I_j:=[1/(j+1),j)$ for
$j=1,2,\dots$. The inner product of $F$ is chosen such that
the functions $e_j$ equal to $1$ over $I_j$ are orthonormal.
Then the algorithm $A_n(f)=\sum_{j=1}^n\il f,e_j\ir_Fe_j$
minimizes the worst case error for all $L_p$ with $p\in[1,\infty)$.
The error is $[n(n+1)]^{-1/p}$.
Since $\il f,e_j\ir_F=f(1/(j+1))$, we may say that this algorithm uses standard
information. Therefore $r^{\rm all-wor}(F,L_p)=
r^{\rm std-wor}(F,L_p)=2/p$.}.
But this is a property of $F$, not the indication of the power of standard
information. By taking the infimum with respect to $F$, we concentrate
on the power of standard information as compared to the power of
linear information.

Suppose now that we take the minimal $n=n^{\,\rm wor-all/std}(\e,F,L_p)$
for which the minimal worst case error
is $\e$ or $\e\,\|I\|$. Assume for simplicity that
$$
e^{\rm all-wor}_n(F,L_p)=n^{-r}\ \ \ \mbox{and}\ \ \
e^{\rm std-wor}_n(F,L_p)=n^{-\alpha}
$$
for some positive $\alpha=r^{\rm std-wor}(F,L_p)\le r$. Then
$$
n^{\,\rm wor-all}(\e,F,L_p)
=\left\lceil\e^{-1/r}\right\rceil
\ \ \ \mbox{and}\ \ \
n^{\,\rm wor-std}(\e,F,L_p)=\left\lceil\e^{-1/\alpha}\right\rceil.
$$
Clearly,
$$
\lim_{\e\to0}\ \frac{\ln\,n^{\,\rm wor-all}(\e,
F,L_p)}{\ln\,n^{\,\rm wor-std}(\e,F,L_p)}
=\frac{\alpha}{r} \ge \ell^{\,\rm wor-x}(r,p).
$$

Hence, if $\ell^{\,\rm wor-x}(r,p)=1$
then function values are as powerful as arbitrary linear functionals.
On the other hand, the smaller $\ell^{\,\rm wor-x}(r,p)$
the less powerful are function values
as compared to arbitrary linear functionals.
If $\ell^{\,\rm wor-x}(r,p)=0$ then the polynomial behavior of
$n^{\,\rm all}(\e,F,L_p)$ in $\e^{-1}$ can be drastically changed
for $n^{\,\rm std}(\e,F,L_p)$.

\begin{rem}
It is well known that, in some cases, we can restrict ourselves
only to linear algorithms. This holds when
$p=\infty$ or when $F$ is a Hilbert space.
Then the corresponding infima for the minimal worst case errors
are attained by
$$
A_n(f)=\sum_{j=1}^nL_j(f)h_j
$$
for some $L_j\in \Lambda\in\{\lstd,\lall\}$ and $h_j\in L_p$.
Much more about the
existence
of linear optimal error algorithms can be found
in e.g.,~{\rm \cite{NW08}}.
\end{rem}

\subsection{Double Hilbert Case}  \label{s2.1}

In this subsection, we consider the approximation problem
defined over a Hilbert
space with the error measured also in the Hilbert space $L_2$.
That is why the name of this subsection is the double Hilbert case.
Approximation in the $L_2$ norm for Hilbert spaces has been studied in
many papers. For our problem the most relevant papers are
\cite{HNV08}, \cite{KWW09} and \cite{WW01}.

Assume that $H$ is a Hilbert space of functions defined on
a set $\Omega$. Since we assume that
function values are continuous this means that $H$ is a reproducing
kernel Hilbert space, $H=H(K)$, where $K$ is defined on
$\Omega\times\Omega$.
Let $L_2=L_2(\Omega,\mu)$ be
the space of $\mu$-square integrable functions with a
measure $\mu$ on $\Omega$. Since the embedding
$I:H(K)\to L_2(\Omega,\mu)$ is continuous, we have
$$
\int_{\Omega}|f(t)|^2\, {\rm d}\mu(t)<\infty\ \ \
\mbox{for all}\ \ \
f\in H(K).
$$
In particular, we can take $f=K(\cdot,t)$ for arbitrary $t\in \Omega$,
since such a function $f$
belongs to $H(K)$. Therefore $W=I^*I:H(K)\to H(K)$,
where $I^*$ is  defined by
$\il g,I(f)\ir_{L_2(\Omega,\mu)}=\il I^*(g),f\ir_{H(K)}$
for all $f \in H(K)$ and $g \in L_2(\Omega, \mu)$,
is given by
$$
W(f)\,(x)=\int_{\Omega}K(x,t)\,f(t)\,{\rm d}\mu(t)\ \ \
\mbox{for all}\ \ \ f\in H(K).
$$
The operator $W$ is self-adjoint and positive semi-definite.
It is well known that
$$
\lim_ne^{\rm wor-all}_n(H,L_2)=0
$$
if and only if $W$ is compact,
see, e.g., \cite[Section 4.2.3]{NW08}.
Unfortunately, in general, $W$ needs not be compact
and therefore  $e^{\rm wor-all}_n(H,L_2)$
does not have to go to zero.
In fact, the sequence $e^{\rm wor-all}_n(H,L_2)$ can be an arbitrary
non-increasing sequence as the following example shows.

\begin{example}[Arbitrary Sequence $e^{\rm wor-all}_n(H,L_2)$]\

Let $(\alpha_n)_{n\in \nat}$ be an arbitrary non-increasing
sequence of nonnegative numbers. Define $k^*$ as the number of
positive $\alpha_n$. If all $\alpha_n$ are positive, we formally set
$k^*=\infty$. If $k^*$ is finite let $\nat_{k^*}=\{1,2,\dots,
k^*\}$,
otherwise let $\nat_{k^*}=\nat$.

For $k\in \nat_{k^*}$,
take arbitrary disjoint intervals $I_k$ of positive
Lebesgue measure $|I_k|$ such that
$\bigcup_{k\in \nat_{k^*}}I_k=[0,1]$, and define the functions
$e_k:[0,1]\to\reals$ by
$$
e_k=\frac{\sqrt{\alpha_k}}{\sqrt{|I_k|}}\ 1_{I_k},
$$
where $1_{I_k}$ is the indicator function of $I_k$.
That is, $e_k(x)=\sqrt{\alpha_k/|I_k|}$ for $x\in I_k$ and $e_k(x)=0$
for $x\notin I_k$.

Define the Hilbert space $H=\overline{{\rm span}}
\{e_k\ |\ k\in\nat_{k^*}\}$ equipped with
the inner product such that $\il e_k,e_j\ir_{H}=\delta_{k,j}$ for all
$k,j\in\nat_{k^*}$.
This means that $H$ is the space of piecewise constant functions
$f:[0,1]\to\reals$ such that
$$
f=\sum_{k=1}^{k^*}a_ke_k\ \ \ \mbox{with}\ \ \ a_k=\il f,e_k\ir_H\ \
\mbox{and}\ \
\|f\|_H=\bigg(\sum_{k=1}^{k^*}a_k^2\bigg)^{1/2}<\infty.
$$
The Hilbert space $H$ has the reproducing kernel
$$
K(x,y)=\sum_{k=1}^{k^*}e_k(x)e_k(y)\ \ \ \mbox{for all}\ \ \
x,y\in[0,1].
$$
Indeed, first of all note that $K$ is well defined since for all $x$
and $y$ the last series has at most one nonzero term.
Then $\il K(\cdot,y_i),K(\cdot,y_j)\ir_{H}=K(y_i,y_j)$, and
$$
0\le\bigg\|\sum_{j=1}^ma_jK(\cdot,y_j)\bigg\|^2_H
=\sum_{i,j=1}^m  a_ia_jK(y_i,y_j).
$$
This shows that the matrix $(K(y_i,y_j))_{i,j=1,2,\dots,m}$ is
symmetric and positive semi-definite for all $m$ and $y_j$. Clearly,
$$
\il f,K(\cdot,y)\ir_{H}=\sum_{k=1}^{k^*}a_ke_k(y)=f(y),
$$
and this completes the proof of the fact
that $K$ is the reproducing kernel of $H$.

Let $L_2=L_2([0,1])$ be the usual space of square Lebesgue integrable
functions. Note that
$$
\|e_k\|_2=\frac{\alpha_k}{\sqrt{|I_k|}}\left(\int_{I_k}\,{\rm
  d}t\right)^{1/2}=\alpha_k.
$$
Therefore,  for any $f\in H$,  we have
$$
\|I(f)\|_2=\|f\|_2=\left(\sum_{k=1}^{k^*} a_k^2\alpha_k^2\right)^{1/2}\le
\alpha_1\|f\|_{H}.
$$
The last bound is sharp, and therefore $\|I\|=\alpha_1$ showing
that $H$ is continuously embedded in $L_2$. The operator $W$ takes now
the form
$$
W(f)=\sum_{k=1}^{k^*}\il f,e_k\ir_2\,e_k.
$$
Note that $W(e_k)=\|e_k\|^2_2\,e_k=\alpha_k^2\,e_k$.
This means that $(\alpha_k^2,e_k)$ are the eigenpairs of $W$
and
$$
W(f)=\sum_{k=1}^{k^*} \alpha_k^2\,\il f,e_k\ir_{H}\,e_k.
$$
It is well known that
$$
e^{\rm wor-all}_n(H,L_2)
=\alpha_{n+1}\ \ \ \mbox{for all} \qquad n \in \nat ,
$$
see, e.g., \cite[Section 4.2.3]{NW08}.
This proves that the behavior of $e^{\rm wor-all}_n(H,L_2)$ can be
arbitrary and, in general, we do not have convergence
of $e^{\rm wor-all}_n(H,L_2)$ to zero.
Clearly, $W$ is compact if and only if $\lim_n\alpha_n=0$.

In addition, this example also shows that for a given $\beta\ge0$
we can define a sequence $\alpha_k$ such that
$r^{\rm all-wor}(H,L_2)=\beta$. Indeed, it is enough to take
$\alpha_k=k^{-\beta}$.  \qed
\end{example}

We discuss the power function $\ell^{\,\rm wor-H}$.
We now assume that $r^{\rm all-wor}(H,L_2)=r>0$.
In particular,  we assume
that the operator $W$ is compact. Then $W$ has eigenpairs
$(\lambda_j,\eta_j)$,
$$
W(\eta_j)=\lambda_j\eta_j\ \ \ \mbox{for all}\ \ j=1,2,\dots,
$$
with $\il \eta_j,\eta_k\ir_{H}=\delta_{j,k}$.
Without loss of generality, we can order the eigenvalues $\lambda_j$
such that $\lambda_1\ge\lambda_2\ge\cdots$. For all $f\in H$, we have
$$
\il f,\eta_k\ir_2=\il I(f),I(\eta_k)\ir_2=\il
f,W(\eta_k)\ir_H=\lambda_k\il f,\eta_k\ir_H.
$$
In particular, letting $f=\eta_j$, we conclude
that the functions $\eta_j$ are also orthogonal in the space $L_2$.

As above,  it is well known that
$$
e^{\rm wor-all}_n(H,L_2)
=\sqrt{\lambda_{n+1}}\ \ \
\mbox{for all}\ \ \ n\in\nat.
$$

If $(e^{\rm wor-all}_n(H,L_2))$ is convergent to zero then
the same also holds for function values, i.e.,
$(e^{\rm wor-std}_n(H,L_2))$
is also convergent to zero.
Indeed, we can reason as in Section 10.4 of~\cite{NW08} that
all linear functionals can be approximated with an arbitrarily
small error when we use function values, and then it is enough to
remember that the error $\sqrt{\lambda_{n+1}}$ is achieved by a linear
algorithm that uses the
$n$ linear functionals $\il f,\eta_j\ir_{H(K)}$.

We have
$$
{\rm trace}(W):=\sum_{j=1}^\infty\lambda_j=
\int_{\Omega}K(x,x)\,{\rm d}\mu(x)=
\sum_{n=0}^\infty
\left[e^{\rm wor-all}_n(H,L_2)\right]^2
$$
and this is finite if $r^{\rm all-wor}(H,L_2)>\tfrac12$.
If $r^{\rm all-wor}(H,L_2)=\tfrac12$ then
$\sum_{n=0}^\infty
\left[e^{\rm wor-all}_n(H,L_2)\right]^2$
may be finite or infinite, and
if $r^{\rm all-wor}(H,L_2)<\tfrac12$ then
$\sum_{n=0}^\infty
\left[e^{\rm wor-all}_n(H,L_2)\right]^2$
is infinite.

The result from~\cite{KWW09} states that
$r^{\rm all-wor}(H,L_2) =r > \tfrac12$
implies
$$
r^{\rm std-wor}(H,L_2) \ge
r-\frac{r}{2r+1}=\frac{2r^2}{2r+1}.
$$

The case
$\sum_{n=0}^\infty
\left[e^{\rm wor-all}_n(H,L_2)\right]^2=\infty$
was studied in~\cite{HNV08}.
It was shown that for any $r\in[0,\tfrac12]$
there is a Hilbert space $H$ such that
$$
r^{\rm all-wor}(H,L_2) = r \quad
\hbox{and} \quad
r^{\rm std-wor} (H, L_2) =0 .
$$

These results give us the following bounds on the power function
$\ell^{\,\rm wor-H}(\cdot,2)$.

\begin{thm}[\cite{HNV08,KWW09}]  \label{thm1}
\begin{eqnarray*}
\ell^{\,\rm wor-H}(r,2)&=&0\ \ \ \ \ \ \ \ \ \ \ \ \
\ \ \ \ \ \mbox{for all}\ \ \
r\in(0,\tfrac12],\\
\ell^{\,\rm wor-H}(r,2)&\in&\left[\frac{2r}{2r+1},1\right]
\ \ \ \ \, \mbox{for all}\ \ \
r\in(\tfrac12,\infty).
\end{eqnarray*}
\end{thm}

Although we do not know the power function $\ell^{\,\rm wor-H}(\cdot,2)$
exactly, we know that there is a jump at $\tfrac12$ since
$\ell^{\,\rm wor-H}(r,2)\ge1/2$ for all $r>1/2$.
Note also that for large $r$, the values of $\ell^{\,\rm wor-H}(r,2)$ are
close to $1$. This means that the power of function values
for $r\in(0,\tfrac12)$ is zero, and is almost optimal for large $r$.

The problem of finding the exact values of $\ell^{\,\rm wor-H}(r,2)$
for $r>\tfrac12$ is one of the main open problems
in the worst case setting.
We know that many people, including the two of us, spent a lot of time
trying to solve this problem but so far in vain. That is why we propose
an open problem with the hope that it
will  soon be solved by the reader.

\begin{OP}
Suppose that $r>\tfrac12$.
Is it true that
$$
\ell^{\,\rm wor-H}(r,2) = 1\, ?
$$
If not, what are  the values of $\ell^{\,\rm wor-H}(r,2)$?
\end{OP}

The rate of convergence neglects to distinguish between
sequences that differ
by a power of logarithms of $n$. Indeed, for $c_n=n^{-r}$
and $b_n=n^{-r}[\ln\,(n+1)]^\beta$ for a positive $r$ and
an arbitrary~$\beta$, we have
$r(c_n)=r(b_n)=r$ independent   of $\beta$. Obviously, for some
standard spaces, we would like to know not only the rate but also
a power of logarithms. We discuss this point in the next example,
where we use the notation
$$
c_n\asymp b_n
$$
which means that there exist
positive numbers $a_1$ and $a_2$ such that $a_1\le c_n/b_n\le a_2$
for large~$n$.

\begin{example}[Sobolev spaces, $p=2$]\

a) For the standard Sobolev spaces $W_2^s([0,1]^d)$ with
an arbitrary $s>0$, which measures the total smoothness of functions,
it is well known that
$$
e^{\rm all-wor}_n(W_2^s([0,1]^d),L_2) \asymp n^{-s/d}.
$$
Of course, in general, function values are not well defined in
$W_2^s([0,1]^d)$. We must assume the embedding condition
$2s>d$ and then function values are well defined and they are
continuous linear functionals. Furthermore, it is known that
$$
e^{\rm all-wor}_n(W_2^s([0,1]^d),L_2) \asymp
e^{\rm std-wor}_n(W_2^s([0,1]^d),L_2) \asymp n^{-s/d},
$$
see,  e.g., {\rm \cite{NW08}} for a survey of such results.

b) For the Sobolev spaces $W_2^{r, {\rm mix}}([0,1]^d)$
with $r >0$, which measures the smoothness of functions
with respect to each variable, it is known that
$$
e^{\rm all-wor}_n(W_2^{r, {\rm mix}}([0,1]^d),L_2)
\asymp n^{-r} (\log n )^{(d-1)r},
$$
see, e.g., \cite{Ga96,
MW81,SU09,Te93,TWW88,Wo94},
where this result can be found in various generalities.

For function values, we must assume that $r>1/2$, and
then the best upper bound is
$$
e^{\rm std-wor}_n(W_2^{r, {\rm mix}}([0,1]^d),L_2)
=\mathcal{O}\left(n^{-r} (\log n)^{(d-1)(r+1/2)}
\right) ,
$$
see {\rm \cite{SU09,Te93,T10}}.

It is \emph{not} known whether this extra power $(d-1)/2$ of logarithms
is needed. It would be very interesting to verify whether
$$
e^{\rm all-wor}_n(W_2^{r, {\rm mix}}([0,1]^d),L_2)\asymp e^{\rm
  std-wor}_n(W_2^{r, {\rm mix}}([0,1]^d),L_2)
$$
holds also for this example. \qed
\end{example}

The examples in \cite{HNV08} use very irregular sequences
$(e_n^{\rm all-wor} (H, L_2))$ and hence do not exclude
a positive answer to the question in the next open problem.

\begin{OP}
Assume that
$e_n^{\rm all-wor} (H, L_2) \asymp n^{-r} \, [ \ln (n+1)]^\beta$
with arbitrary $r >0$ and $\beta \in \reals$.
Is it true that this implies
$$
e_n^{\rm std-wor} (H, L_2) \asymp
e_n^{\rm all-wor} (H, L_2) ?
$$
\end{OP}

\subsection{Single Hilbert Case}

In this short subsection, we mostly consider the approximation problem
defined over a Hilbert space with the error measured in the
non-Hilbert space $L_p$ for $p\not=2$. That is why the name of this
subsection is the single Hilbert case.

We report on a recent result of Tandetzky~\cite{Ta10} who
considered the approximation problem for arbitrary $p\in[1,\infty)$.
He proved that
for any $r\in(0,\min(\tfrac1p,\tfrac12)]$ there exists a Hilbert space
$H$ continuously embedded in $L_p=L_p([0,1])$ such that
$$
r^{\rm all-wor}(H,L_p) =r  \qquad \hbox{and} \qquad
r^{\rm std-wor}(H, L_p) = 0 .
$$
This result obviously implies that the power function
is zero over $(0,\min(\tfrac1p,\tfrac12)]$.
It seems to us that no example is known in the literature
for a Hilbert space for which
$e^{\rm all-wor}_n(H,L_p)$  tends to zero faster than the
sequence $e^{\rm std-wor}(H,L_p)$
with the additional assumption that
$r^{\rm all-wor}(H,L_p)  >\min(\tfrac1p,\tfrac12)$.
This implies that we do not know the behavior of the power function
over $(\min(\tfrac1p,\tfrac12),\infty)$. We summarize our partial
knowledge of the power function in the following theorem.

\begin{thm}[\cite{Ta10}]    \label{thm2}
Let $p\not=2$.
\begin{eqnarray*}
\ell^{\,\rm wor-H}(r,p)&=&0\ \ \ \ \ \ \ \ \mbox{for all}\ \ \
r\in(0,\min(\tfrac1p,\tfrac12)] .  
\end{eqnarray*}
\end{thm}

Only for the case $p=\infty$ do we know a little more
about the behavior of the power function.
In this case the rates are related as explained in the following theorem.

\begin{thm}[\cite{No88}]   \label{thm3}
Let $F$ be a Hilbert or a Banach space.
Then
\begin{equation}   \label{eq1}
e^{\rm std-wor}_n(F, L_\infty ) \le (1+n) \,
e^{\rm all-wor}_n (F, L_\infty)\ \ \ \mbox{for all}
\ \ \ n\in\nat .
\end{equation}
\end{thm}

This inequality follows from Proposition 1.2.5, page 16,
in {\rm \cite{No88}}, where it is stated for the
Kolmogorov widths and also applies to the linear or Gelfand widths.

The inequality \eqref{eq1} cannot be improved even if we assume
that $F$ is a Hilbert space. This follows from the following example.

\begin{example}   \label{ex3}
Take $F=H=\reals^{n+1}$. That is, $f\in H$ is now defined on
$\{1,2,\dots,n+1\}$ and can be identified with $f=[f_1,f_2,
\dots,f_{n+1}]$, where $f_i=f(i)$.
The space $H$ is equipped with the inner product
$$
\il f,g\ir_{H}=\left[\sum_{i=1}^{n+1}f_i\right]
\left[\sum_{i=1}^{n+1}g_i\right]
+\eps\,\sum_{i=1}^{n+1}f_ig_i\ \ \
\mbox{for all}\ \ \ f,g\in H.
$$
The unit ball of $H$ is thus
$$
\textstyle{
B= \bigl\{ f \in \reals^{n+1} \mid \ \
\left[\sum_{i=1}^{n+1} f_i\right]^2 + \eps
\sum_{i=1}^{n+1} f_i^2 \le 1 \bigr\}
}.
$$
Then for $\eps \to 0$, we obtain
$$
e^{\rm std-wor}_n (F, L_\infty ) \ge1.
$$
Indeed, knowing $f(x_i)$ for $i=1,2,\dots,n$, with
$x_i\in\{1,2,\dots,n+1\}$, we take $f$ such that $f(x_i)=0$.
Since we have at most $n$ conditions on $n+1$ components of $f$
then at least one component of $f$ from the unit ball
is free and can be taken as
$\pm 1/\sqrt{1+\eps}$.
This proves that the worst case error of any
algorithm is at least $1/\sqrt{1+\e}$ which in the limit as
$\eps$ goes to zero is $1$.

Consider the information
$$
N(f) = \left[f_1-f_2, f_2-f_3, \dots , f_{n} - f_{n+1}
 \right]
\ \ \ \mbox{for all}\ \ \ f\in H.
$$
It is known that the minimal error of all algorithms that use $N$
is the supremum of $\|f\|_H$ for $f\in B$ and $N(f)=0$.
Observe that $N(f)=0$ implies that
$f= [c,c, \dots , c]$. Next, $f\in B$ implies that
$$
c^2 \le \frac{1+\eps/(n+1)}{(n+1)^2}.
$$
Hence,  again for $\eps \to 0$,
we obtain $e^{\rm all-wor}_n (F, L_\infty )\le1/(n+1)$. \qed
\end{example}

Let $r^{\rm all-wor} (F,L_\infty) = r >1$.
Then the
inequality~\eqref{eq1} implies that
$$
r^{\rm std-wor}(F,L_\infty) \ge r-1 .
$$
Thus, Theorem~\ref{thm3} implies the following
behavior of the power function for $p=\infty$.

\begin{thm}   \label{thm4}
\begin{eqnarray*}
\ell^{\,\rm wor-H/B}(r,\infty)&\in& \left[\frac{r-1}r,1\right]
\ \ \ \ \ \mbox{for all}
\ \ \ r>1.
\end{eqnarray*}
\end{thm}

Hence, for both $p=2$ and $p=\infty$,
we see that for large $r$, the power function is almost one.

We want to guess the behavior of the power function for
$r>\min(\tfrac1p,\tfrac12)$.
It can be helpful to see the actual rates of convergence for
some standard spaces. In particular, for $p=\infty$, the rates
are known for Sobolev spaces.

\begin{example}[Sobolev spaces, $p=\infty$] \

a) For the Sobolev spaces $W_2^s([0,1]^d)$ and an arbitrary $s$ for
which $2s>d$,
it is well known that
$$
e^{\rm all-wor}_n(W_2^s([0,1]^d),L_\infty) \asymp
e^{\rm std-wor}_n(W_2^s([0,1]^d),L_\infty) \asymp
n^{-s/d+1/2} ,
$$
see, e.g., {\rm \cite{NW08}}.

\goodbreak

b) For the Sobolev spaces $W_2^{s, {\rm mix}}([0,1]^d)$
with $s > 1/2$, it is known that
$$
e^{\rm all-wor}_n(W_2^{s, {\rm mix}}([0,1]^d),L_\infty) \asymp
e^{\rm std-wor}_n(W_2^{s, {\rm mix}}([0,1]^d),L_\infty) \asymp
n^{-s+1/2} (\log n )^{(d-1)s},
$$
see {\rm \cite{T93}}. \qed
\end{example}

Hence, at least for the standard Sobolev spaces the rates
are the same even up to logarithmic factors.
This again suggests that the power function can be just one
for all $r > (\min(\tfrac1p,\tfrac12),\infty)$.
This is the next open problem.

\begin{OP}
Verify whether it is true that for all $p\in[1,\infty]$ we have
$$
\ell^{\,\rm wor-H}(r,p)=\left\{
\begin{array}{ll}
0 &\mbox{for all $r\in\big(0,\min(\tfrac1p,\tfrac12)\big]$,}\\
1 &\mbox{for all $r\in\big(\min(\tfrac1p,\tfrac12),\infty\big)$.}
\end{array}\right.
$$
\end{OP}

We end this section with a remark on the rates of convergence for
different $p$.

\begin{rem}
It is interesting to compare the sequences
$$
e^{\rm all-wor}_n(H, L_p) \qquad \mbox{and/or} \qquad
e^{\rm std-wor}_n(H,L_p)
$$
for the same $H$ but different $p$.
The following example shows that, in general, there exists \emph{no}
relation between these sequences.
Some relations do exist as shown in \cite{KWW08}
but under some additional assumptions about $H$.
The following example shows that
some assumptions on $H$ are indeed needed,
otherwise everything can happen.

Take $L_2=L_2([0,1])$, $L_\infty = L_\infty([0,1])$ and assume
that $[0,1]$ is the disjoint union of intervals
$I_k$ of positive length $\lambda_k$
such that
$\sum_{k=1}^\infty \lambda_k =1$. Assume also that
$$
\lambda_1 \ge \lambda_2 \ge \cdots
$$
and put $e_k = 1_{I_k}$.
We define a Hilbert space $H$ by its unit ball
$$
B=\biggl\{ \sum_{k=1}^\infty \alpha_k e_k  \
\bigg| \  \sum_{k=1}^\infty \frac{\alpha_k^2}{\gamma_k^2}
\le 1 \biggr\} ,
$$
where
$$
\gamma_1 \ge \gamma_2 \ge \dots > 0
\ \ \ \mbox{with}\ \ \
\lim_{k \to \infty} \gamma_k = 0.
$$
Hence for
$f = \sum_{k=1}^\infty \alpha_k e_k \in H$, we obtain
$$
\Vert f \Vert_H^2 = \sum_{k=1}^\infty
\frac{\alpha_k^2}{\gamma_k^2}  \qquad
\hbox{and} \qquad
\Vert f \Vert_2^2  = \sum_{k=1}^\infty \alpha_k^2 \, \lambda_k,\ \ \
\ \ \Vert f \Vert_\infty = \sup_k |\alpha_k|.
$$
{}From this, we easily conclude that the optimal approximation
for $L_2$ as well as for $L_\infty$ is given by
$$
f = \sum_{k=1}^\infty \alpha_k e_k\  \mapsto\
\sum_{k=1}^n \alpha_k e_k .
$$
Note that
$$
\alpha_k=\il f,e_k\ir_H=f(x_k)\,\lambda_k,
$$
where $x_k\in I_k$. This means that the optimal error algorithm
for function values and linear functionals is the same, and therefore
$$
e^{\rm all-wor}(H,L_p)=e^{\rm std-wor}(H,L_p)\ \ \ \mbox{for}
\ \ \ p\in\{2,\infty\}.
$$
However,
$$
e^{\rm all-wor}_n
(H, L_\infty ) = \gamma_{n+1}
\qquad \hbox{and} \qquad
e^{\rm all-wor}_n (H, L_2) = \gamma_{n+1}\,\sqrt{\lambda_{n+1}}.
$$
Since $\{\gamma_n\}$ and $\{\lambda_n\}$ are not related, it
is easy to get an example with
$$
r^{\rm all-wor}(H,L_\infty) =0 \qquad \hbox{but} \qquad
r^{\rm all-wor}(H,L_2) = \infty .
$$
Hence, in general, the difference between the minimal rates
for $L_2$ and $L_\infty$ approximation can be extreme.
\end{rem}

\subsection{Banach Case}  \label{s2.2}

In this subsection, we study the approximation problem
defined over a Banach space that is continuously embedded in $L_p$.
As always, we assume that function evaluations
are continuous functionals.
We establish some bounds on the power functions by recalling known
results for Sobolev spaces.

\begin{example}[Sobolev spaces, $1 \le p<\infty$]\

For the Sobolev space $W_p^s([0,1]^d)$ for an arbitrary $s>0$,
it is known that
$$
e^{\rm all-wor}_n(W_p^s([0,1]^d),L_p)\asymp n^{-s/d}.
$$
Function values are well defined in $W_p^s([0,1]^d)$ only
if the embedding condition $s/d>1/p$
or $s=d$ and $p=1$ holds.
However,
we may use the approach suggested in~\cite{He06c} that
allows us to consider the case without this embedding condition.
Namely, we
limit ourselves only to continuous functions by taking
$$
F=W_p^s([0,1]^d)\cap C([0,1]^d)
$$
with norm
$$
\|f\|_F=\|f\|_{W_p^s([0,1]^d)}+\|f\|_{C([0,1]^d)}.
$$
Here, $C([0,1]^d)$ is the space of continuous functions equipped
with the max norm. Then $F$ is a Banach space for which
function values are well defined and
function values are continuous linear functionals
on this space. Then
for $s/d\le 1/p$
and $s/d < 1$ in the case $p=1$,
respectively,  it was shown in {\rm \cite{He06c}} that
$$
e^{\rm std-wor}_n (F,L_p)\asymp 1.
$$
\vskip-\baselineskip\qed
\end{example}

The last example implies that
\begin{eqnarray*}
\ell^{\,\rm wor-B}(r,p)&=&0 \qquad \mbox{for all} \qquad
r\in(0,1/p] \quad \mbox{and} \quad 1 < p < \infty ,  \\
\ell^{\,\rm wor-B}(r,1)&=&0
\qquad   \mbox{for all} \qquad r \in (0,1).
\end{eqnarray*}

We now show that $\ell^{\,\rm wor-B}(r,p)=0$
over larger domains of $r$ for a given $p$
by recalling other results for Sobolev spaces.

\begin{example}[Sobolev space {$W_1^s([0,1]^d)$}, $1
\le p < \infty$] \label{ex6} \

Consider
the approximation problem for the Sobolev space
$W^s_{1}([0,1]^d)$ with error measured in $L_p=L_p([0,1]^d)$.
This problem is well defined
and convergent for the class $\lall$
if we assume that
$s/d> 1-1/p$.

For $p\in[1,2]$, we have
$$
e^{\rm all-wor}_n(W_{1}^s([0,1]^d),L_p)\asymp n^{-s/d},
$$
whereas for $p\in[2,\infty)$, we have
$$
e^{\rm all-wor}_n(W_{1}^s([0,1]^d),L_p)\asymp n^{-s/d+1/2-1/p},
$$
see e.g.,~\cite{Vy08}.
The last relation also holds for $p=\infty$ as will be
needed later.

The same results are also valid for the space
$F=W^s_{1}([0,1]^d)\cap C([0,1]^d)$
with the norm
$$
\|f\|_F=\|f\|_{W_{1}^s([0,1]^d)}+\|f\|_{C([0,1]^d)}.
$$
For the space $F$, we can consider function values
for all $s/d>1-1/p$. For $s/d\le1$, we have
$$
e^{\rm std-wor}_n(F,L_p)\asymp 1.
$$\vskip-\baselineskip\qed
\end{example}

Let $p\in[1,2]$.
The previous example implies that
$$
\ell^{\,\rm wor-B}(r,p)=0 \qquad \mbox{for all} \qquad
r\in\left(1-\frac1p, \ 1\right].
$$
For $p \in (1,2]$, we showed before
that $\ell^{\,\rm wor-B}(r,p)=0$ for all
$r\in(0,1/p]$. Since $(0,1/p]\cup(1-1/p,1]=(0,1]$, we
obtain
$$
\ell^{\,\rm wor-B}(r,p)=0 \qquad  \mbox{for all} \qquad
r\in(0,1] \quad \mbox{and} \quad   p\in[1,2].
$$

Let $p\in[2,\infty)$. The previous example implies that
$$
\ell^{\,\rm wor-B}(r,p)=0 \qquad \mbox{for all} \qquad
r\in\left(\frac12,\  \frac12+\frac1p\right].
$$
Now we show
that $\ell^{\, \rm wor-B}(r,p)=0$ also for $p \in [2,
\infty )$ and $r \in (0,\tfrac{1}{2}]$.
We increase the space
$F=W^s_{1}([0,1])\cap C([0,1])$
with the norm
$$
\|f\|_F=\|f\|_{W_{1}^s([0,1])}+\|f\|_{C([0,1])}
$$
(for $d=1$) even more by adding functions
from a H\"older class $C^\alpha $, where
$0 < \alpha \le 1/2$.
Hence we take the space
$$
\wt F = F + C^\alpha
$$
with the norm
$$
\Vert f \Vert_{\wt F} :=
\inf \{ \Vert g \Vert_F + \Vert h \Vert_{C^\alpha} \mid \
f= g+h, \
g \in F, \
h \in C^\alpha \} .
$$
Since the unit ball of $\wt F$ is larger than that of $F$
we still have
$e_n^{\rm std-wor} (\wt F, L_p) \asymp 1$ for $s \le 1$.
It is well known that
$e_n^{\rm all-wor} (C^\alpha ([0,1]) , L_p) \asymp n^{-\alpha}$
and the same holds for $\wt F$ if $\alpha \le  s-1/2+1/p$.
Hence for $p \in [2, \infty)$, we obtain
$$
\ell^{\,\rm wor-B}(r,p)=0 \qquad \mbox{for all} \qquad
r\in\left(0  ,\  \frac12+\frac1p\right].
$$

We learnt some properties of the power function by using
known results for Sobolev spaces $W^s_{p_1} ([0,1]^d)$
in the case  $s/d\le 1/p_1$
so that function values did not even supply
convergence. Since we needed to assume that $s/d>1/p_1-1/p$,
the case $p=\infty$ could not be covered.

We now recall some results for Sobolev spaces
when the embedding condition is satisfied and when there is a
difference in the convergence rates between function values and
arbitrary linear functionals.

\begin{example}[Sobolev space {$W_1^s([0,1]^d)$, $1 \le p \le \infty$}] \

Consider the approximation problem for the Sobolev space $W_1^s([0,1]^d)$
with error measured in $L_p$.
We now assume that $s/d\ge1$. Then function values are well defined
and are continuous linear functionals. Furthermore,
$$
e^{\rm std-wor}_n (W_1^s([0,1]^d),L_p)\asymp n^{-s/d+1-1/p},
$$
see,  e.g.,  the survey of such results in Section 4.2.4 of~\cite{NW08}
or
\cite{Vy07,Vy08}.
\qed
\end{example}

The last two examples  imply  the following estimates of the power
function. For all $r>1$ and $p\in[1,2]$, we have
$$
\ell^{\,\rm wor-B}(r,p)\le 1-\frac1r\left(1-\frac1p\right),
$$
and for all $r>1$ and $p\in[2,\infty]$, we have
$$
\ell^{\,\rm wor-B}(r,p)\le 1-\frac1{2r}.
$$

We summarize the properties of the power function established
in this section in the following theorem.
The only case where we have a positive lower bound
is the case $p=\infty$,  see Theorem~\ref{thm4}.

\begin{thm}    \label{thm5}
\begin{eqnarray*}
\ell^{\,\rm wor-B}(r,p)&=&0\ \ \ \ \ \ \ \ \
\ \ \ \ \ \ \ \qquad \quad
\mbox{for all}\ \ \
r\in(0,1]\ \ \mbox{and}\ \ p\in [1, 2]    ,\\
\ell^{\,\rm wor-B}(r,p)&=&0\ \ \ \ \ \ \ \ \
\ \ \ \ \ \ \ \qquad \quad
\mbox{for all}\ \ \
r\in(0, \tfrac{1}{2} + \tfrac{1}{p}]  \ \
\mbox{and}\ \ p\in(2, \infty )     ,\\
\ell^{\,\rm wor-B}(r,p)&\le& 1-\frac1r\left(1-\frac1p\right)
\ \ \ \ \ \ \, \mbox{for all}\ \ \ r>1\ \ \mbox{and}\ \ p\in[1,2],\\
\ell^{\,\rm wor-B}(r,p)&\le& 1-\frac1{2r}\ \ \ \ \ \ \ \ \ \ \ \ \
\ \ \ \ \ \mbox{for all}\ \ \ r>1\ \ \mbox{and}\ \ p\in[2,\infty), \\
1- \frac{1}{r} \le \ell^{\,\rm wor-B}(r,\infty)&\le& 1-\frac1{2r}\ \ \
\ \
\ \ \ \ \ \ \  \ \  \ \ \  \ \mbox{for all}\ \ \ r>1 .
\end{eqnarray*}
\end{thm}

It is interesting to note that although we do not know the exact
values of the power functions in the Hilbert and Banach cases,
we can check that they are different at least for $p=2$.
Indeed, from Theorems~1 and~5, we have
\begin{eqnarray*}
\ell^{\,\rm wor-B}(r,2)&=&
\ell^{\,\rm wor-H}(r,2)\qquad\qquad\ \ \
\ \ \ \mbox{for all}\ \ \  r\in(0,\tfrac12],\\
\ell^{\,\rm wor-B}(r,2)=
0&<&\tfrac12\le \ell^{\,\rm wor-H}(r,2)\quad\qquad\ \ \
\mbox{for all}\ \ \ r\in(\tfrac12,1],\\
\ell^{\,\rm wor-B}(r,2)\le1-\frac1{2r}&<&\frac{2r}{2r+1}\le
\ell^{\,\rm wor-H}(r,2)\ \ \ \ \ \mbox{for all}\ \ \
r\in(1,\infty).
\end{eqnarray*}
This shows that at least for $p=2$ the power of function values
for the Hilbert case is larger than for the Banach case for all
$r>\tfrac12$.

Obviously, it would be desirable to find the exact values
of the power function $\ell^{\,\rm wor-B}(r,p)$
for all $r\in(0,\infty)$ and $p\in[1,\infty]$.
However, it could be a very difficult problem.
Hence, as maybe a less difficult problem,
we would like to check
the following property of the power function.
\begin{OP}
For $p\in[1,\infty]$, find the supremum $a^*(p)$ of $a$ for which
$$
\ell^{\,\rm wor-B}(r,p)=0\ \ \ \mbox{for all}\ \ \ r\in(0,a].
$$
We only know that $a^*(p)\ge1$ for all $p\in[1,\infty)$.
\end{OP}

We already indicated that the power functions for the Hilbert and
Banach cases are different for $p=2$. It would be of interest to
check if this holds for all $p$.
\begin{OP}
Find all $p\in[1,\infty]$ for which
$$
\ell^{\,\rm wor-B}(\cdot,p)\not=
\ell^{\,\rm wor-H}(\cdot,p).
$$
\end{OP}
\vskip 1pc
Similar to  Example~\ref{ex3}, we present an example of
a Banach space $F$ where the ratio
$$
\frac{e^{\rm std-wor}_n (F, L_p )}{e^{\rm all-wor}_n (F, L_p )}
$$
is large for $p>1$ and a fixed $n$.

\begin{example}   \label{ex9}
Take $F=\ell^{n+1}_1$, i.e.,
$F=\reals^{n+1}$ with the $\ell_1$ norm. Then we obtain
\begin{equation}    \label{eq99}
e^{\rm std-wor}_n (F, L_p ) = (n+1)^{1-1/p}
e^{\rm all-wor}_n (F, L_p ) ,
\end{equation}
since
$e^{\rm std-wor}_n (F, L_p ) = 1$ and
$e^{\rm all-wor}_n (F, L_p ) = (n+1)^{1/p-1}$.
The upper bound in the last statement follows again
with the information
$N(x)= (x_2-x_1, x_3-x_2, \dots , x_{n+1}-x_n)$
while the lower bound follows from the fact that the unit ball
of $\ell^{n+1}_1$ contains a $\ell_p^{n+1}$ ball of radius
$(n+1)^{1/p-1}$.

Again this ratio $(n+1)^{1-1/p}$ as in \eqref{eq99}
can be obtained with a Hilbert space and actually we can take the same
spaces as in Example~\ref{ex3}, i.e., we define in $H=\reals^{n+1}$
the scalar product
$$
\il f,g\ir_{H}=\left[\sum_{i=1}^{n+1}f_i\right]
\left[\sum_{i=1}^{n+1}g_i\right]
+\eps\,\sum_{i=1}^{n+1}f_ig_i\ \ \
\mbox{for all}\ \ \ f,g\in H,
$$
and consider the limit where $\eps >0$ tends to zero.   \qed
\end{example}

We end this section with another open problem.

\begin{OP}
Find the supremum of
$e^{\rm std-wor}_n (F, L_p ) /
e^{\rm all-wor}_n (F, L_p )$
over all Banach and/or Hilbert spaces.
So far, we know that
\begin{equation}
\sup_F  \frac{e^{\rm std-wor}_n (F, L_p )}
{e^{\rm all-wor}_n (F, L_p )}  \ge (n+1)^{1-1/p} ,
\end{equation}
and equality holds if $p=\infty$.
\end{OP}

\section{Randomized setting}  \label{s3}

We approximate the embedding operator $I: F \to L_p$ in the
randomized setting. We now briefly define this setting.
The reader may find more
on this  subject, e.g., in~\cite{NW08,NW10,TWW88}.

We approximate $I$ by algorithms $A_{n}$ that use
$n$ values of linear functionals on the average
and each linear functional
is chosen randomly
with respect to a probability distribution.

More precisely,
the algorithm $A_{n}$ is of the following form
\begin{equation}\label{157}
A_{n}(f,\omega)=\phi_{n,\omega}\left(
L_{1,\omega_{1}}(f),
L_{2,\omega_{2}}(f),
\dots,
L_{n(\omega),\omega_{n(\omega)}}(f)
\right),
\end{equation}
   and  the number $n(\omega)$ of functionals can also be random.
Here $\omega=[\omega_{1},\omega_{2},\dots]$, and
the linear functionals
$L_{j,\omega_{j}}$ are random functionals distributed according to
a probability distribution on elements
$\omega_{j}$ which may depend on $j$ as well as on the
values already computed, i.e., on $L_{i,\omega_i}(f)$ for
$i=1,2,\dots, j-1$.
The mapping $\phi_{n,\omega}:\R^{n(\omega)}\to L_p$
is a random mapping, and
$$
\expect_{\omega}\,n(\omega)\le n.
$$
We also allow adaptive choices of the functionals
$L_{j, \omega_j}$.  That is,
$L_{j,\omega_j}$ may depend on the already selected functionals
and the values $L_{1, \omega_1} (f),
L_{2, \omega_2} (f) , \dots ,
L_{j-1, \omega_{j-1}}(f)$.

Without loss of generality, we assume that $A_{n}(f,\cdot)$ is
measurable, and define the randomized error of $A_{n}$ as
$$
e^{\rm ran}(A_{n})=\sup_{\|f\|_F \le 1 }
\left(\expect_{\omega} \Vert I(f)-A_{n}(f,\omega) \Vert_p^2\right)^{1/2}.
$$

Again, we compare such algorithms with algorithms that are based
on function values, i.e.,
each
$L_{j, \omega_j}$ is now of the form
$L_{j, \omega_j}(f) = f(t_{j, \omega_j})$
and
\begin{equation}\label{158}
A_{n}(f,\omega)=\phi_{n,\omega}\left(f(t_{1,\omega_{1}}),f(t_{2,
\omega_{2}}),
\dots,f(t_{n(\omega),\omega_{n(\omega)}})\right).
\end{equation}

Hence, we consider algorithms that use $n$ linear functionals either
from the class $\lstd$ or the class $\lall$.
We define the minimal errors as follows.

\begin{defi}
For $n\in \nat_0$, let
$$
\enrall(F, L_p) =
\inf\left\{e^{\rm ran}(A_{n})\,|\ \
L_j \in \lall  \ \mbox{and}\ A_{n}\  \mbox{as in~\eqref{157}} \right\},
$$
and
$$
\enrstd(F, L_p) =
\inf\left\{e^{\rm ran}(A_{n})\,|\ \
 L_j \in \lstd\ \mbox{and}\
A_{n}\  \mbox{as in~\eqref{158}} \right\}.
$$
\end{defi}

As in the worst case setting, for $n=0$ it is easy to see that the best
algorithm is $A_{0}=0$ and obtain
$$
e^{\rm all-ran}_0(F,L_p)= e^{\rm std-ran}_0(F, L_p)=
\sup_{\Vert f \Vert_F \le 1} \Vert f \Vert_p =
\sup_{\Vert f \Vert_F \le 1} \Vert I(f) \Vert_p =
\Vert I \Vert.
$$
This is the initial error that can be achieved without
computing any linear functional on the functions $f$.
Clearly,
$$
\enrall(F, L_p ) \le \enrstd(F, L_p)\ \ \ \mbox{for all} \ \ \  n\in\nat.
$$
The sequences $\left(\enrall(F, L_p)\right)$
and $\left(\enrstd(F, L_p)\right)$
are both non-increasing but not necessarily convergent to zero.

As in the worst case setting,
we want to compare the rates of convergence
$$
r^{\rm all-ran} (F, L_p)
= r\left(\enrall(F,L_p)\right) \qquad
\hbox{and}  \qquad
r^{\rm std-ran} (F, L_p)
= r\left(\enrstd(F,L_p)\right) .
$$
In particular, we would like to know if it is
possible that the sequence
$\left(r^{\rm all-ran} (F, L_p)\right)$ converges
much faster than the sequence
$\left(r^{\rm std-ran} (F, L_p)\right)$.
The main question addressed in this section is to find or estimate
the~\emph{power} function
defined as
$\ell^{\rm \, ran-x}:(0,\infty)\times[1,\infty]\to[0,1]$ by
$$
\ell^{\,\rm ran-x}(r,p)
:=\inf_{F:\, r^{\rm all-ran}(F,L_p)=r
}\frac{r^{\rm std-ran}(F,L_p)}r,
$$
where ${\rm x}\in\{H,B\}$     indicates that
the infimum is taken over all Hilbert spaces $({\rm x}=H)$ or over all
Banach spaces $({\rm x}=B)$ continuously embedded in $L_p$
and the rate of convergence is $r$ when we use
arbitrary linear functionals.
In the randomized setting, we do \emph{not}
need to assume
that function values are continuous linear functionals.

\subsection{Double Hilbert Case}  \label{s3.1}

In this subsection, we consider the approximation problem
defined over a Hilbert
space with the error measured also in the Hilbert space $L_2$.
It may be surprising but
the results in the double Hilbert case
are complete due to \cite{WW07},
and there is no need to discuss different cases
depending on the values of $r$.

\begin{thm}[\cite{WW07}]\label{thm6}
Let
$ I: H \to L_2 (\Omega)$
be a continuous
embedding from a Hilbert space $H$ into
$L_2(\Omega)$. Then
$$
r^{\rm all-ran}(H,L_2) =
r^{\rm std-ran}(H,L_2) .
$$
Therefore
$$
\ell^{\rm \, ran-H} (r,2) = 1 \qquad
\hbox{for all} \quad r > 0.
$$
\end{thm}
\vskip 1pc
We add that it was known before, see \cite{No92,Was89}, that also
$$
r^{\rm all-ran}(H,L_2) =
r^{\rm all-wor}(H,L_2) .
$$
This means that the power of function values in the randomized setting
is the same as the power of arbitrary linear functionals in the worst
case setting, which in turn is the same as in the randomized setting.

\subsection{Other Cases}  \label{s3.2}

For $p>2$, we know examples from the literature where the rate
$
r^{\rm all-ran}(H,L_p)$
is larger
than the rate
$
r^{\rm std-ran}(H,L_p)
$.
Namely take
$I: W^r_2([0,1]) \to L_p ([0,1])$.
Then with $\lall$ one can achieve the order $n^{-r}$
(with additional log terms in the case $p=\infty$,
but the order is still $r$), see \cite{Ma91}.
For
$\lstd$ the optimal order is $n^{-r+1/2-1/p}$, see \cite{He06c}.
The authors of~\cite{He06c,Ma91}
studied the case of integer $r$,
but the results can be extended
via interpolation to all $r > 1$. Therefore, we obtain
$$
\ell^{\rm \, ran-H} (r,p) \le \frac{r-1/2+1/p}{r} \qquad
\hbox{if} \quad r \ge 1 \quad \hbox{and} \quad  p>2.
$$
We summarize these estimates of the power
function in the following theorem.

\begin{thm}
Let $p>2$. Then
\begin{eqnarray*}
\ell^{\rm \, ran-H}(r,p)&\le&1-\frac{1/2-1/p}{r}\ \ \ \ \
\mbox{for all} \ \ \ r \ge 1 .
\end{eqnarray*}
\end{thm}

Sobolev embeddings in the randomized setting
were studied by several authors,  including
\cite{He06c,He09a,He09b,Ma91,No88,TWW88,Was89}.
For our purpose,
the most important papers are \cite{He06c,Ma91} and
the paper~\cite{He09a} for the interpolation argument.

For the embedding
$I: W^r_2([0,1]) \to L_\infty ([0,1])$
the rate is improved by 1/2 if we switch
from the class $\Lambda^{\rm std}$  to the class $\Lambda^{\rm all}$.
This gap of 1/2 is the largest possible
under some additional conditions, see \cite{KWW08,KWW09b}.
Let us add in passing that
the same gap of 1/2
appears for $\lall$ between the worst case and the randomized setting.

The Hilbert case for $p\in[1,2)$ as well as the Banach case
for all $p\in[1,\infty]$ have not yet been studied. We pose
this as an open problem.

\begin{OP}
Study the power function in the randomized setting for the Hilbert
case with $p\in[1,2)$ and for the Banach case for all
$p\in[1,\infty]$. In particular, determine the supremum $a^*(p)$
of $a$ for which
$$
\ell^{\rm \,ran-H/B}(r,p)=0\ \ \ \mbox{for all}\ \ \ r\in(0,a].
$$
\end{OP}

\section{Average case setting with a Gaussian measure}  \label{s4}

In the average case setting, we assume
that $I : F \to L_p (\Omega)$ is continuously embedded and
function evaluations are continuous functionals on $F$.
As far as we know, only the case $p=2$ was studied and we
report the known results from \cite{HWW08}  for this case.

We assume that $F$ is a separable Hilbert/Banach space equipped with
a zero mean Gaussian measure $\mu$.
As in the worst case setting, we consider deterministic algorithms,
and due to general results, see \cite{TWW88}, it is enough
to compare linear algorithms
$$
A_n(f) = \sum_{k=1}^n L_k(f) g_k  \quad \hbox{and} \quad
A_n(f) = \sum_{k=1}^n f(x_k) g_k ,
$$
where $g_k \in L_2(\Omega)$.
The average case error of an algorithm is defined by
$$
e^{\rm avg} (A) := \biggl(  \int_F
\Vert f - A(f) \Vert_p^2 \, {\rm d} \mu (f) \biggr)^{1/p} .
$$
As in the other settings, we define the minimal $n$th average case errors
$
\enaall(F, L_p)
$,
$\enastd(F, L_p)$
and the power function $\ell^{\rm \,avg-H/B}$. That is,
for
$$
r^{\rm all/std-avg}(F,L_p)=r(e_n^{\rm all/std-avg}(F,L_p))
$$
we have
$$
\ell^{\rm \,avg-x}(r,p):=\inf_{F:\,r^{\rm all-avg}(F,L_p)=r}
\frac{r^{\rm std-avg}(F,L_p)}r.
$$
As always, ${\rm x}\in\{H,B\}$  and we take the infimum over separable
Hilbert (${\rm x}=H$)
or Banach (${\rm x}=B$) spaces equipped with zero mean Gaussian
measures that are continuously embedded in $L_p$ and for which
function values are continuous linear functionals as well as the rate of
convergence is $r$ when arbitrary linear functionals are used.

As already mentioned, results are known only for $p=2$.
Then the cases of the Hilbert and Banach spaces are  the
same due to the presence of Gaussian measures. This follows from the
fact that even if $F$ is a separable Banach space then
the minimal errors for the class $\lall$
depend on the Gaussian measure $\nu=\mu\,I^{-1}$ given by
$$
\nu(M)=\mu\left(\left\{f\in F\,|\ \ I(f)\in M\right\}\right\}
$$
for a Borel set $M$ of $L_2$. The measure $\nu$ is also a zero mean
Gaussian measure whose covariance operator $C_\nu:L_2\to L_2$
is given by
$$
\il C_\nu f_1,f_2\ir_{L_2}=\int_{L_2}\il f,f_1\ir_{L_2}\il
f,f_2\ir_{L_2}
\, {\rm d}\nu(f)\ \ \ \mbox{for all}\ \ \ f_1,f_2\in L_2.
$$
The operator $C_\nu$ is self adjoint, positive semi-definite, compact
and has a finite trace. That is, its ordered eigenvalues $\lambda_j$
have a finite sum. It is known that
$$
e_n^{\rm
  all-avg}(F,L_2)=\bigg(\sum_{j=n+1}^\infty\lambda_j\bigg)^{1/2}.
$$

As in the randomized setting
for the double Hilbert space,
the results on the power function are complete
and there is no need to discuss
different cases of $r$.

\begin{thm}[\cite{HWW08}]\label{thm8}
Let
$ I: F \to L_2 (\Omega)$
be a continuous embedding from a separable Banach space $F$ equipped with
a zero mean  Gaussian measure $\mu$ into  $L_2(\Omega)$.
Then
$$
r^{\rm all-avg}(F,L_2) =
r^{\rm std-avg}(F,L_2) .
$$
Therefore
$$
\ell^{\rm \, avg-H/B}  (r,2) = 1 \qquad
\hbox{for all} \quad r > 0.
$$
\end{thm}
Of course it would be interesting to study the power function
for other values of $p$. This is posed as our last open problem.

\begin{OP}
Study the power function in the average case setting for $p\not=2$.
In particular, verify whether a similar result as Theorem 8 holds.
\end{OP}

\vskip 1pc
\noindent
{\bf Acknowledgment}\

We appreciate comments on this paper from
Stefan Heinrich,  Anargyros Papageorgiou, Joseph F. Traub, Grzegorz
W. Wasilkowski and two anonymous referees.
We  especially thank  Stefan Heinrich for
pointing out some errors in our previous manuscript.

E.N.
was partially supported by
the DFG-Priority Program 1324.
H.W. was partially
supported by the National Science
Foundation.




{

\bigskip

\hskip1.4 em\vbox{\noindent
Erich Novak \\
Mathematisches Institut, Universit\"at Jena \\
Ernst-Abbe-Platz 2, 07740 Jena, Germany \\
{\tt novak@mathematik.uni-jena.de} \\
{\tt http://users.minet.uni-jena/$\sim$novak}}

\medskip

\hskip1.4 em\vbox{\noindent
Henryk Wo\'zniakowski \\
Department of Computer Science, Columbia University\\
New York, NY 10027, USA, and\\
Institute of Applied Mathematics, University of Warsaw\\
ul. Banacha 2, 02-097 Warszawa, Poland\\
{\tt henryk@cs.columbia.edu} \\
{\tt http://www.cs.columbia.edu/$\sim$henryk}}

}

\endddoc
\begin{thebibliography}{99}

\setlength{\parsep }{-0.5ex}
\setlength{\itemsep}{-0.5ex}

\frenchspacing

\newcommand\BAMS{\emph{Bull. Amer. Math. Soc.\ }}
\newcommand\BIT{\emph{BIT\ }}
\newcommand\Com{\emph{Computing\ }}
\newcommand\CA{\emph{Constr. Approx.\ }}
\newcommand\FCM{\emph{Found. Comput. Math.\ }}
\newcommand\JAT{\emph{J. Approx. Th.\ }}
\newcommand\JC{\emph{J. Complexity\ }}
\newcommand\JMA{\emph{SIAM J. Math. Anal.\ }}
\newcommand\JMAA{\emph{J. Math. Anal. Appl.\ }}
\newcommand\JMM{\emph{J. Math. Mech.\ }}
\newcommand\MC{\emph{Math. Comp.\ }}
\newcommand\NM{\emph{Numer. Math.\ }}
\newcommand\RMJ{\emph{Rocky Mt. J. Math.\ }}
\newcommand\SJNA{\emph{SIAM J. Numer. Anal.\ }}
\newcommand\SR{\emph{SIAM Rev.\ }}
\newcommand\TAMS{\emph{Trans. Amer. Math. Soc.\ }}
\newcommand\TOMS{\emph{ACM Trans. Math. Software\ }}
\newcommand\USSR{\emph{USSR Comput. Maths. Math. Phys.\ }}

\frenchspacing


\bibitem{Ga96}
E. M. Galeev,
Linear widths of H\"older-Nikolskii classes
of periodic functions of several variables,
Math. Notes {\bf 59}, 133--146, 1996.

\bibitem{He06c}
S. Heinrich,
Randomized approximation of Sobolev embeddings,
in:  \emph{Monte Carlo and Quasi-Monte Carlo Methods 2006},
A. Keller, S. Heinrich, H. Niederreiter (eds.),
445--459,
Springer, Berlin, 2008.

\bibitem{He09a}
S. Heinrich,
Randomized approximation of Sobolev embeddings II,
\JC  {\bf 25}, 455--472, 2009.

\bibitem{He09b}
S. Heinrich,
Randomized approximation of Sobolev embeddings III,
\JC  {\bf 25}, 473--507, 2009.

\bibitem{HWW08}
F. Hickernell, G. W. Wasilkowski and
H. Wo\'zniakowski,
Tractability of linear multivariate problems in the
average case setting,
in:  \emph{Monte Carlo and Quasi-Monte Carlo Methods 2006},
A. Keller, S. Heinrich, H. Niederreiter (eds.),
461--494,
Springer, Berlin, 2008.

\bibitem{HNV08}
A. Hinrichs, E. Novak and J. Vyb\'\i ral,
Linear information versus function evaluations
for $L_2$-approximation,
\JAT {\bf 153}, 97--107, 2008.

\bibitem{KWW08}
F. Y. Kuo, G. W. Wasilkowski and H. Wo\'zniakowski,
Multivariate $L_\infty$
approximation in the worst case setting
over reproducing kernel Hilbert spaces,
\JAT  {\bf 152}, 135--160, 2008.

\bibitem{KWW09}
F. Y. Kuo, G. W. Wasilkowski and H. Wo\'zniakowski,
On the power of standard information for multivariate
approximation in the worst case setting,
\JAT  {\bf 158}, 97--125, 2009.

\bibitem{KWW09b}
F. Y. Kuo, G. W. Wasilkowski and H. Wo\'zniakowski,
On the power of standard information for $L_\infty$
approximation in the randomized case setting,
\emph{BIT Numer. Math.}
{\bf 49}, 543--564, 2009.

\bibitem{Ma91}
P. Math\'e,
Random approximation of Sobolev embeddings,
\JC {\bf 7}, 261--281, 1991.

\bibitem{MW81}
C. A. Micchelli and G. Wahba,
Design problems for optimal surface interpolation,
in
\emph{Approximation Theory and Applications},
Z. Ziegler ed., pp. 329--347,
Academic Press, New York, 1981.

\bibitem{No88}
E. Novak,
\emph{Deterministic and Stochastic Error Bounds in
Numerical Analysis},
LNiM {\bf 1349}, Springer-Verlag, Berlin, 1988.

\bibitem{No92}
E. Novak,
Optimal linear randomized methods for linear operators
in Hilbert spaces,
\JC {\bf 8}, 22--36, 1992.

\bibitem{NW08}
E. Novak and H. Wo\'zniakowski,
\emph{Tractability of Multivariate Problems},
Volume I: Linear Information, European Math. Soc., Z\"urich,
2008.

\bibitem{NW10}
E. Novak and H. Wo\'zniakowski,
\emph{Tractability of Multivariate Problems},
Volume II: Standard Information for Functionals,
EMS, Z\"urich, 2010, to appear.

\bibitem{Pie80}
A. Pietsch,
\emph{Operator Ideals},
North Holland, 1980.

\bibitem{SU09}
W. Sickel and T. Ullrich,
Spline interpolation on sparse grids,
to appear in Applicable Analysis.

\bibitem{Ta10}
R. Tandetzky,
Approximation of functions
from a Hilbert space using function values or
general linear information, in progress, 2010.

\bibitem{Te93}
V. N. Temlyakov,
\emph{Approximation of Periodic Functions},
Nova Science,
New York, 1993.

\bibitem{T93}
V. N. Temlyakov,
On approximate recovery of functions with bounded mixed derivative.
\JC {\bf 9}, 41--59, 1993.

\bibitem{TWW88}
J. F. Traub, G. W. Wasilkowski and H. Wo\'zniakowski,
\emph{Information-Based Complexity},
Academic Press, 1988.

\bibitem{T10}
H. Triebel,
\emph{Bases in Function Spaces, Sampling, Discrepancy, Numerical
Integration},
EMS Publ. House, Z\"urich, 2010.

\bibitem{Vy07}
J. Vyb\'{\i}ral,
Sampling numbers and function spaces,
\JC {\bf 23}, 773--792, 2008.

\bibitem{Vy08}
J. Vyb\'{\i}ral,
Widths of embeddings in function spaces,
\JC {\bf 24}, 545--570, 2008.

\bibitem{Was89}
G. W. Wasilkowski,
Randomization for continuous problems,
\JC {\bf 5}, 195--218, 1989.

\bibitem{Wo94}
H. Wo\'zniakowski,
Tractability and strong tractability of multivariate tensor product
problems,
\emph{J. of Computing and Information} {\bf 4}, 1--19, 1994.

\bibitem{WW01}
G. W. Wasilkowski and H. Wo\'zniakowski,
On the power of standard information for weighted
approximation,
\emph{Found. Comput. Math.}
{\bf 1}, 417--434, 2001.

\bibitem{WW07}
G. W. Wasilkowski and H. Wo\'zniakowski,
The power of standard information for multivariate
approximation in the randomized case setting,
\emph{Math. Comp.}
{\bf 76}, 965--988, 2007.

\end{thebibliography}
